\documentclass[11pt]{article}

\usepackage{amsfonts,amsmath,latexsym}
\usepackage{amstext,amssymb}
\usepackage{amsthm}
\usepackage{color}
\usepackage{soul}

\newcommand {\R}{\mathbb{R}}

\newcommand {\B}{\mathbf{B}}

\newcommand {\calCp}{\mathcal{C}_p}

\newtheorem {thm} {Theorem}

\newtheorem {lemma}{Lemma}

\newcommand{\beq}{\begin{equation}}
\newcommand{\eeq}{\end{equation}}
\numberwithin{equation}{section}

\begin {document}

\title{ 
{\bf A reverse H\"older inequality for extremal Sobolev functions}}
\author{Tom Carroll\footnote{School of Mathematical Sciences, University 
College Cork, {\tt t.carroll@ucc.ie}} 
\ and Jesse Ratzkin\footnote{Department of Mathematics and Applied 
Mathematics, University of Cape Town, {\tt jesse.ratzkin@uct.ac.za}}
}
\maketitle

\begin {abstract}\noindent
Let $n \geq 2$, let $\Omega \subset \R^n$ be a bounded domain with
$\mathcal{C}^1$ boundary, and let $1 \leq p < \frac{2n}{n-2}$ (simply $p\geq 1$ if 
$n=2$). The well-known Sobolev imbedding theorem and Rellich compactness implies 
$$\calCp(\Omega) = \inf \left \{ \frac{\int_\Omega |\nabla f|^2 dm}{\left (
\int_\Omega |f|^p dm \right )^{2/p},  } : f \in W^{1,2}_0(\Omega) , f \not \equiv 0\right \}$$
is a finite, positive number, and the infimum is achieved by a nontrivial extremal function 
$u$, which one can assume is positive inside $\Omega$. We prove that for every $q > p$ 
there exists $K= K(n, p, q, \calCp(\Omega))>0$ such that $\| u \|_{L^p(\Omega)} \geq K \| u \|_{L^q(\Omega)}.$
This inequality, which reverses the classical H\"older inequality, mirrors results of 
G. Chiti \cite {Chiti2} for the first Dirichlet eigenfunction of the Laplacian and of M.\ van den Berg 
\cite{vdB} for the torsion function. 

\end {abstract}

\section{Introduction and statement of results}
In 1972, Payne and Rayner \cite{PR} showed that the eigenfunction $\phi$ of the Dirichlet Laplacian 
corresponding to the first eigenvalue $\lambda(D)$ for a bounded planar domain $D$ 
satisfies a reverse H\"older inequality, specifically 
\[
\int_D \phi^2\,dm \leq\frac{\lambda(D)}{4\pi}  \left( \int_D \phi\,dm\right)^2,
\]
with $dm$ being Lebesgue measure. 
The inequality is isoperimetric in the sense that equality holds if and only if $D$ is a disk. 

It proved not to be entirely straightforward to extend this inequality to regions in higher dimensions. 
Payne and Rayner \cite{PR2} obtained an isoperimetric extension to higher dimensions that they themselves 
described as not \lq entirely satisfactory\rq, since their inequality became trivial for regions 
of given volume but large eigenvalue. 
Kohler-Jobin \cite{KJ} obtained an isoperimetric comparison between  the 
$L^2$ and the $L^1$ norms of the eigenfunction that did not suffer from the 
defects of that of Payne and Rayner. Her inequality is 
\[
\int_\Omega \phi^2\,dm \leq\frac{\lambda(\Omega)^{n/2}}{2n\,\omega_n\, j_{n/2-1}^{n-2}}  
		\left( \int_\Omega \phi\,dm\right)^2,
\]
only leaving unanswered whether equality could hold for regions other than balls. 
Here $\Omega$ is a bounded domain in $\R^n$,  the volume of the unit ball in $\R^n$ is denoted by $\omega_n$, 
and $j_m$ denotes the first positive zero of the Bessel function $J_m$. She 
obtained further results in \cite{KJ3}.
Subsequently, Chiti \cite{Chiti1, Chiti2}  obtained, for $0 < p < q$,  an inequality of the form  
\begin{equation}\label{tc1}
\left( \int_\Omega \phi^q\,dm \right) ^{1/q} \leq  K\big(p,q, n\big) \, \lambda(\Omega)^{\frac{n}{2}(\frac{1}{p}- \frac{1}{q})} 
	\left( \int_\Omega \phi^p\,dm \right) ^{1/p}
\end{equation}
where $\Omega$ is a bounded region in $\R^n$. 
The inequality is isoperimetric in that equality holds if and only if $\Omega$ is a ball. 

The torsion function $u$ of an open set $\Omega$ in $\R^n$ is the solution of the problem 
$\Delta u = -1$ in $\Omega$ with zero Dirichlet boundary conditions. 
For $1 \leq p < q < \infty$, van den Berg \cite[Theorem~1]{vdB} obtains, as part of his work, 
an upper bound for the $L^q(\Omega)$-norm of the torsion function $u$ in terms of its $L^p(\Omega)$-norm, 
by which he concludes that if $u$ is in $L^p(\Omega)$ then it belongs to $L^q(\Omega)$ for 
all $q>p$. (A precursor of this result is Corollary 2 of \cite{vdBC}.) 

Our aim here is to obtain a reverse H\"older inequality for the extremal Sobolev function 
that is similar to that of Chiti for the eigenfunction of the Laplacian.  
For a bounded domain $\Omega$ with $\mathcal{C}^1$ boundary and, for admissible values of $p$, 
we consider 
\begin {equation} \label {sobolev-functional} 
\begin{split}
\calCp (\Omega) & =  \inf \left \{ \Phi_p(f) : f \in W^{1,2}_0(\Omega), 
f \not \equiv 0\right \} \\
& =  \inf \left \{ \frac{\int_\Omega |\nabla f|^2\, dm}
{\left ( \int_\Omega |f|^p\, dm \right )^{2/p} } : f \in W^{1,2}_0(\Omega), 
f \not \equiv 0\right \}.
\end{split}
\end {equation}
Here the allowable range of exponents is $1 \leq p < \frac{2n}{n-2}$ if $n\geq 3$
and $p \geq 1$ if $n=2$. We can see, using the change of variables $y = x/r$, that 
$\calCp$ obeys the scaling law 
\begin {equation} \label{sobolev-scaling1} 
\calCp(r\Omega) = r^{n-2-\frac{2n}{p}} \calCp(\Omega) = r^{\alpha_{n,p}} \calCp(\Omega),
\end {equation}
and observe that within the allowable range of exponents $\alpha_{n,p}< 0$. 
By Rellich compactness and the Sobolev embedding theorem, $\calCp(\Omega)$ 
is a finite, positive number, and is realized by a nontrivial function $u = u_p \in 
W^{1,2}_0(\Omega)$, which we can take to be positive inside $\Omega$. 
This extremal Sobolev function satisfies 
\begin{equation}\label{Sobolev-pde}
\Delta u + \Lambda u^{p-1} = 0, \quad \left. u \right |_{\partial 
\Omega} = 0, \quad \Lambda = \calCp(\Omega) 
\left ( \int_\Omega u^p dm \right )^{\frac{2-p}{p}}.
\end{equation}
Note that $\Phi_2(f)$ is the Rayleigh-Ritz quotient for $f$ so that 
$\mathcal{C}_2(\Omega) = \lambda(\Omega)$, and the extremal 
Sobolev function $u$ is the first Dirichlet eigenfunction $\phi$ of the Laplacian.
The torsional rigidity $P(\Omega)$ of the region $\Omega$ is $4/\mathcal{C}_1(\Omega)$
and is given by $\int_\Omega u_1\,dm$ where $\Delta u = -2$ in $\Omega$ and 
$u$ vanishes on the boundary of $\Omega$. 

Previously, we have generalized the Payne-Rayner theorems to extremal 
Sobolev functions in two dimensions \cite{CRpPR1} and in 
higher dimensions \cite{CRpPR2}. However, our proof suffers the same 
flaw as the one Payne and 
Rayner employed, and our result in \cite{CRpPR2} is trivial for domains 
such that $\calCp$ is large compared to the volume. 

\begin {thm} 
Let $1 \leq p \leq 2$ and let $q  \geq p$. Then there is a positive constant $K = 
K(n, p, q, \calCp(\Omega))$ such that 
\begin {equation} \label {main-ineq} 
\left ( \int_\Omega u^p dm \right )^{1/p} \geq K \left ( \int_\Omega
u^q dm \right )^{1/q}. \end {equation} 
Equality can only occur if $\Omega$ is a round ball, except for a set of measure zero. 
\end {thm}
The restriction $1 \leq p \leq 2$ is necessary for our proof, but we expect a 
similar inequality to hold for all admissible values of $p$. 
The constant $K(n, p,q,\calCp(\Omega))$ is something one can in principle calculate, but 
it is not given by a closed formula. However, we can extract exactly how
$K$ depends on $\calCp(\Omega)$ as follows. Let $\B^*$ be the ball with 
$\calCp(\Omega) = \calCp(\B^*)$ and let $\phi$ be the extremal function for $\calCp(\B^*)$. 
We will see below that 
\begin {equation} \label{optimal-constant1} 
K(n,p,q,\calCp(\Omega)) = \frac{\left (\int_{\B^*} |\phi|^p dm \right )^{1/p}}
{\left ( \int_{\B^*} |\phi|^q dm \right )^{1/q}}.\end {equation}
If we denote the radius of $\B^*$ by $\rho$, then the change of variables 
$\psi(r) = \phi(\rho r)$ gives us an extremal function on the unit ball $\B$, and 
we obtain 
\begin{eqnarray} \label {sobolev-scaling2}
\frac{\left ( \int_{\B^*} |\phi|^p dm \right )^{1/p}}{\left ( \int_{\B^*} 
|\phi|^q dm \right )^{1/q}} & = & \frac{\rho^{n/p} \left (\int_\B
|\psi|^p dm \right)^{1/p}}{\rho^{n/q} \left (\int_\B 
|\psi|^q dm \right )^{1/q}} \\ \nonumber 
& = & \rho^{n\left ( \frac{1}{p} - \frac{1}{q}\right )}
\frac{\left ( \int_\B |\psi|^p dm \right )^{1/p}}{\left ( \int_\B |\psi|^q 
dm \right )^{1/q}} . 
\end {eqnarray}
Using \eqref{sobolev-scaling1} we find 
\begin {equation} \label{sobolev-scaling3}
\rho = \left ( \frac{\calCp(\B^*)}{\calCp(\B)} \right )^{1/\alpha_{n,p}}
= \left ( \frac{\calCp(\Omega)}{\calCp(\B)} \right)^{1/\alpha_{n,p}},
\end {equation} 
and so combining \eqref{sobolev-scaling2} and \eqref{sobolev-scaling3}
yields 
\begin {eqnarray} \label{optimal-constant2} 
K & = & (\calCp(\Omega))^{\frac{n}{\alpha_{n,p}} \left ( \frac{1}{p}- \frac{1}{q} 
\right )}
(\calCp(\B))^{-\frac{n}{\alpha_{n,p}} \left ( \frac{1}{p} - \frac{1}{q} \right )} \frac{\left (
\int_\B |\psi|^p dm \right )^{1/p}}{\left ( \int_\B |\psi|^q dm 
\right )^{1/q}} \\ \nonumber
& = & \hat K(n,p,q) (\calCp(\Omega))^{\frac{n}{\alpha_{n,p}}
\left ( \frac{1}{p}- \frac{1}{q} \right )} 
\end {eqnarray} 
where $\B$ is the unit ball and $\psi$ is an extremal function for $\calCp(\B)$. 
Observe that our dependence on $\mathcal{C}_2(\Omega)$ in \eqref{optimal-constant2}
is the same as the dependence on $\lambda$ in Chiti's inequality \eqref{tc1}, as 
it must be. 
In the special case $p=1$ we see that \eqref{optimal-constant2} 
simplifies to 
$$K(n,1,q,P(\Omega)) = (P(\Omega))^{\frac{n}{n+2} \left ( 1- \frac{1}{q}
\right )} \hat K(n,1,q),$$
where $P(\Omega)$ denotes the torsional rigidity of $\Omega$.

Our technique is inspired by the rearrangements of Talenti \cite{Tal} 
and of Chiti \cite{Chiti2, Chiti1}. However, a key difference is the Chiti 
only considers a PDE which is linear and homogeneous. This allows one to 
scale solutions arbitrarily, choosing either the $L^q$-norm of $u$ to be $1$ for 
some arbitrary $q$, or the $\sup$-norm of $u$ to be $1$. In our case, due to the 
nonlinearity and the inhomogeneity, the 
only natural normalization is the one we choose above, namely $\int u^p dm = 1$. 
With a careful analysis we find that parts of Chiti's proof carries over, and part 
does not. This nonlinearity and inhomogeneity is the main reason we only obtain 
estimates for integrals of certain powers of $u$, rather than the full range of exponents 
as Chiti does in \cite{Chiti1}.

\bigskip \noindent {\sc Acknowledgements:} Much of this research was completed 
while the second author visited the first author at the University College Cork. We 
thank UCC its hospitality. J.~R. is partially 
supported by the National Research Foundation of South Africa. 

\section{Rearrangements} 

In this section we review properties of the distribution function $\mu$ 
associated to our extremal function $u$, and two different symmterizations 
of $u$. 

Associated to $u_p$ (and, indeed, any positive function on $\Omega$) is 
the distribution function, which measures the volume of sublevel sets. Let 
$M = \sup_{x \in \Omega} (u(x))$, and for $0 \leq t \leq M$ define 
\begin {equation} \label {dist-function} 
\mu(t) = |\{ x \in \Omega: u(x)>t\}|.\end {equation} 
By Sard's theorem, for almost every $t$ the set $\Omega_t = \{ u>t\}$ has a smooth 
boundary 
$$\partial \Omega_t = \{ x \in \Omega: u(x) = t\},$$
and hence, in particular, $\mu$ is defined almost everywhere. We also 
observe that a.e. $t$ we have 
\begin {equation} \label {co-area}
\mu(t) = \int_t^M \int_{\partial \Omega_\tau} \frac{d\sigma}{|\nabla u|}
d\tau \Rightarrow \frac {d\mu}{dt} = -\int_{\partial \Omega_t} \frac{d\sigma}
{|\nabla u|} < 0.\end {equation} 
In particular, $\mu$ is non-increasing, so it has a (left) inverse function 
$u^*(s)$ defined by 
\begin {equation} \label{symm-dec-rearrange1} 
u^*(s) = \inf \{ t \geq 0 : \mu (t) < s \} . \end {equation} 
It will later be useful to notice that $u^*(\mu(t)) = t$ implies 
\begin {equation} \label{dist-inverse-fcn}
\frac{du^*}{ds} = \frac{1}{d\mu/dt}, \end {equation} 
which again holds almost everywhere. 

We can now define the symmetric decreasing rearrangement $u^\dagger$ 
of $u$. First let $\Omega^*$ be the ball centered at $0$ with $|\Omega| 
= |\Omega^*|$, and then define 
\begin {equation} \label{symm-dec-rearrange2} 
u^\dagger(x) = u^*(\omega_n |x|^n), \end {equation}
where $\omega_n$ is the volume of a unit ball in $\R^n$. 
As we did with $\Omega$, we define 
$$\Omega^*_t = \{ x \in \Omega^* : u^\dagger(x) > t\}$$ 
and see that, by construction, $|\Omega_t| = |\Omega^*_t|$ for almost 
every $t$. In particular, $u$ and $u^\dagger$ (and $u^*$) have exactly 
the same function values.  Moreover, these two functions are are 
equimeasurable, so for any $q>0$ 
and $t \in [0, M_u]$ we have 
\begin {equation} \label{equimeasurable} 
\int_{u>t} |u|^q dm = \int_{u^\dagger> t} |u^\dagger|^q dm 
= \int_0^{\mu(t)} (u^* (\tau))^q d\tau.\end {equation} 

The symmetric decreasing rearrangement compares functions on $\Omega$ 
to functions on the ball $\Omega^*$. We can also compare to functions defined 
on the ball $\B^*$, which has the same value of $\calCp$ as $\Omega$ 
({\it i.e.} $\calCp(\Omega) = \calCp(\B^*)$). We let $\phi = \phi_p$ be the 
extremal function for \eqref{sobolev-functional} on $\B^*$, normalized 
so that $\| \phi \|_{L^p(\B^*)} = 1$. Then 
\begin {equation} \label {sobolev-pde2} 
\Delta \phi + \calCp(\B^*) \phi^{p-1} = \Delta \phi + \calCp(\Omega) 
\phi^{p-1} = 0, \qquad \left. \phi \right |_{\partial \B^*}  = 0. \end {equation} 
A theorem of Gidas, Ni, and Nirenberg \cite{GNN} tells us $\phi$ is 
radial and $\frac{d\phi}{dr}< 0$, so in particular 
\begin {equation} \label{max-phi} 
\phi(0) = \sup\{ \phi(x) : x \in \B^* \}.\end {equation}  
Also, by the Faber-Krahn inequality, 
$\B^* \subset \Omega^*$, with equality if and only if $\Omega = 
\Omega^* = \B^*$.

\section{Prelimary analysis} 

We first rewrite the extremal function $\phi$ with respect to volume. For 
$0 \leq s \leq |\B^*|$ define $\phi^*(s)$ by 
$$\phi (x) = \phi^*(\omega_n |x|^n).$$
\begin {lemma} The function $\phi^*(s)$ defined by $\phi(x) = \phi^*\
(\omega_n|x|^n)$ satisfies 
\begin {equation} \label{integro-differential1}
(\phi^*)'(s) = -\calCp(\Omega)n^{-2}\omega_n^{-2/n} s^{-2+\frac{2}{n}}
\int_0^s (\phi^*)^{p-1} (t) dt. \end {equation} 
\end {lemma} 

\begin {proof} Our change of variables is given by 
$$s =\omega_n r^n \Rightarrow \frac{d\phi}{dr} = \frac{ds}{dr} 
\frac{d\phi^*}{ds} = 
n \omega_n r^{n-1} (\phi^*)' = n \omega_n \left ( \frac{s}{\omega_n} 
\right )^{\frac{n-1}{n}} (\phi^*)'.$$
The extremal function $\phi$ satisfies 
$$-\calCp(\Omega) \phi^{p-1} = \Delta \phi = r^{1-n} \frac{d}{dr} \left (
r^{n-1} \frac{d\phi}{dr} \right ),$$
which we can rewrite as 
\begin {eqnarray*} 
-\calCp(\Omega) (\phi^*)^{p-1} & = & \left ( \frac{s}{\omega_n} \right )^{\frac{1-n}{n}} 
n \omega_n \left ( \frac{s}{\omega_n} \right )^{\frac{n-1}{n}} \frac{d}{ds} 
\left ( \left ( \frac{s}{\omega_n} \right )^{\frac{n-1}{n}} n \omega_n 
\left ( \frac{s}{\omega_n} \right )^{\frac{n-1}{n}} (\phi^*)' \right ) \\ 
& = & n^2 \omega_n^{2/n} \frac{d}{ds} \left ( s^{2-\frac{2}{n}} (\phi^*)' 
\right ). 
\end {eqnarray*} 
Rearranging this last equation yields 
$$-\frac{d}{ds} \left ( s^{2-\frac{2}{n}} (\phi^*)' \right ) = \frac{\calCp(\Omega)}
{n^2 \omega_n^{2/n}} (\phi^*)^{p-1}, $$ 
which we can integrate to get 
$$-s^{2-\frac{2}{n}} (\phi^*)'(s) = \frac{\calCp(\Omega)}{n^2 \omega_n^{2/n}}
\int_0^s (\phi^*)^{p-1}(t) dt.$$
The lemma follows.
\end {proof}

The following is essentially equation (34) of \cite{Tal}, but we 
include its proof for the reader's convenience. 
\begin {lemma} 
The function $u^*$ satisfies 
\begin {equation} \label {integro-differential2} 
-(u^*)'(s) \leq \calCp(\Omega) n^{-2} \omega_n^{-2/n} 
s^{-2 + \frac{2}{n}} \int_0^s (u^* (t))^{p-1} dt. \end {equation} 
\end {lemma} 

\begin {proof} 
First use Gauss's divergence theorem to see that for almost every $t$ we have 
\begin {equation} \label {ibp1} 
\calCp(\Omega) \int_{\Omega_t} u^{p-1} dm = -\int_{\Omega_t} \Delta u 
dm = \int_{\partial \Omega_t} \left \langle  \nabla u, 
\frac{\nabla u}{|\nabla u|} \right \rangle d\sigma= \int_{\partial \Omega_t}
|\nabla u| d\sigma.\end {equation}
Combining \eqref{ibp1}, \eqref{co-area}, and the Cauchy-Schwarz inequality, 
we see 
$$|\partial \Omega_t| = \int_{\partial \Omega_t} d\sigma \leq \left [ 
- \mu'(t) \int_{\partial \Omega_t} |\nabla u| d\sigma \right ]^{1/2}
= \left [- \calCp(\Omega) \mu'(t) \int_{\Omega_t} u^{p-1} dm \right ]^{1/2}.$$
Now square this last inequality and apply the isoperimetric inequality to 
obtain 
\begin {equation} \label{integro-differential3}
- \calCp(\Omega) \mu'(t) \int_{\Omega_t} u^{p-1} dm \geq |\partial 
\Omega_t|^2 \geq 
n^2 \omega_n^{2/n} |\Omega_t|^{\frac{2n-2}{n}} = n^2 \omega_n^{2/n} 
(\mu(t))^{2-\frac{2}{n}}.\end {equation}
Finally we make the change of variables $s = \mu(t)$ 
and recall 
$$\frac{d\mu}{dt} = \frac{1}{du^*/ds}, \qquad \int_{\Omega_t} u^{p-1} dm
= \int_0^s (u^*(\tau))^{p-1} d\tau, $$
which turns \eqref{integro-differential3} into 
$$- \frac{1}{(u^*)'} \calCp(\Omega) \int_0^s (u^*(t))^{p-1} dt \geq n^2 
\omega_n^{2/n} s^{2-\frac{2}{n}}.$$
 \end {proof} 

Finally, we quote a result due to Hardy, Littlewood, and P\'olya \cite{HLP}. 
\begin {lemma} \label{HLP-lemma}
Let $M>0$, let $0 < q_1 \leq q_2$, and let 
$f, g \in L^{q_2}([0,M])$. Also let $f^*$ be the decreasing rearrangement 
of $f$ and $g^*$ the decreasing rearrangement of $g$. If 
$$\int_0^s (f^*(t))^{q_1} dt \leq \int_0^s (g^* (t))^{q_1} dt \textrm{ for all }
0 \leq s \leq M$$ 
then 
$$\int_0^M |f|^{q_2} ds \leq \int_0^M |g|^{q_2} ds.$$
\end {lemma}

\section {Proof of the main theorem}

In this section we prove \eqref{main-ineq}. 

If $|\Omega| = |\B^*|$ then, by the Faber-Krahn inequality, we must have 
$\Omega = \Omega^* = \B^*$. In this case, we take 
$$K(n,p,q, \calCp(\Omega)) =  \frac{\left ( \int_{\B^*} \phi^p dm \right )^{1/p}  }
{\left ( \int_{\B^*} \phi^q dm \right )^{1/q} }= 
\left ( \int_{\B^*} \phi^q dm \right )^{-1/q}.$$ 
This will, in fact, turn out to be the optimal constant in the general case. 

We may now assume $\B^*$ is strictly contained in $\Omega^*$. By our 
normalization, 
\begin {equation} \label {integral-ineq1} 
1 = \int_{\B^*} \phi^p dm = \int_0^{|\B^*|} (\phi^*)^p ds = 
\int_0^{|\Omega|} (u^*)^p ds > \int_0^{|\B^*|} (u^*)^p ds, \end {equation} 
so we cannot have $\phi^* \leq u^*$ on the entire interval $[0, |\B^*|]$. 
On the other hand, $\phi^*(|\B^*|) = 0 < u^*(|\B^*|)$. Combining these 
two facts implies that the graphs of $\phi^*$ and $u^*$ must 
cross. Define 
$$s_1 = \inf \{ s \in [0, |\B^*|]: \phi^*(\bar s) < u^*(\bar s) \textrm{ for 
all } \bar s \in (s, |\B^*|] \}, $$ 
and observe that (by continuity) $\phi^*(s_1) = u^*(s_1)$, 
while $\phi^* < u^*$ on $(s_1, |\B^*|]$. The point $s_1$ is the 
first crossing point of the two graphs, as viewed from the right hand side. 

Our next task is to prove that $s_1$ is in fact the only crossing. It is 
useful to first observe $s_1>0$. Indeed, if $s_1 = 0$ then $\int_0^{|\B^*|}
(u^*)^p ds \geq \int_0^{|\B^*|} (\phi^*)^p ds$, which contradicts 
\eqref{integral-ineq1}.

\begin {lemma} On the interval $[0, s_1]$ we have $\phi^*(s) \geq 
u^*(s)$.  \end {lemma} 

\begin {proof} Suppose there exists $s^* \in [0,s_1)$ such that 
$u^*(s^*) > \phi^*(s^*)$, and define 
$$s_2 = \sup \{ s < s_1 : u^* (s) > \phi^*(s) \}.$$
We're supposing that the set $\{s \in [0,s_1) : u^*(s) > \phi^*(s) \}$
is nonempty, so the continuity of $\phi^*$ and $u^*$ implies $s_2>0$. 
Next let 
$$s_3 = \inf \{ s < s_2 : u^*(\bar s) > \phi^*(\bar s) \textrm{ for all }
\bar s \in (s, s_2)\}.$$

Now we have $\phi^* < u^*$ in $(s_1, |\B^*|]$, $\phi^* \geq u^*$ 
in $[s_2, s_1]$, and $\phi^* < u^*$ in $(s_3, s_2)$. Define the function 
$$v^*(s) = \left \{ \begin {array}{rl} \max \{ \phi^*(s), u^*(s) \}  & 0 \leq s \leq s_3 \\ 
u^*(s) & s_3 \leq s \leq s_2 \\ \phi^*(s) & s_2 \leq s \leq |\B^*|. 
\end {array} \right. $$
By construction, we have $v^*> 0$ on $[0, |\B^*|)$ and 
$v^*(|\B^*|) = 0$. We claim  
\begin {equation} \label {integro-differential6} 
-(v^*)' (s) \leq \calCp(\Omega) n^{-2} \omega_n^{-2/n} 
s^{-2 + \frac{2}{n}} \int_0^s (v^*(t))^{p-1} dt\end {equation}
as well. 

For $0 \leq s \leq s_3$ we have $v^*(s) = \max \{ \phi^*(s), 
u^*(s)\}$. This is the maximum of two nonincreasing functions, so 
it is nonincreasing itself, and therefore $v^*$ is differentiable almost 
everywhere. If $u^*(s) \leq \phi^*(s)$ then 
\begin {eqnarray*} 
-(v^*)'(s) &  =&  -(\phi^*)'(s) = \calCp(\Omega) n^{-2} 
\omega_n^{-2/n} s^{-2 + \frac{2}{n}} \int_0^s (\phi^*(t))^{p-1} 
dt \\
& \leq &\calCp(\Omega) n^{-2} \omega_n^{-2/n} s^{-2 + \frac{2}{n}}
\int_0^s (v^*(t))^{p-1} dt, \end {eqnarray*}
and, similiarly, if $u^*(s) > \phi^*(s)$ then 
\begin {eqnarray*}
-(v^*)'(s) & = & -(u^*)'(s) \leq \calCp(\Omega) n^{-2} 
\omega_n^{-2/n} s^{-2 + \frac{2}{n}} \int_0^s (u^*(t))^{p-1} 
dt \\ 
& \leq & \calCp(\Omega) n^{-2} \omega_n^{-2/n} s^{-2 + \frac{2}{n}}
\int_0^s (v^*(t))^{p-1} dt.\end {eqnarray*} 
We conclude \eqref{integro-differential6} holds for $0 \leq s < s_3$. 
Nearly identical arguments show \eqref{integro-differential6} also holds 
in the intervals $(s_3, s_2)$ and $(s_2, |\B^*|)$. 

Associated to $v^*$ we have the radial test function $v: \B^* \rightarrow 
\R$ defined by $v(x) = v^*(\omega_n |x|^n)$. Observe that 
\begin {equation} \label {test-fcn4} 
\int_{\B^*} v^p dm = \int_0^{|\B^*|} (v^*(t))^p dt > \int_0^{|\B^*|} 
(\phi^*(t))^p dt = 1\end {equation} 
(the inequality follows from $v^* > \phi^*$ on $(s_3, s_2)$) and 
\begin {eqnarray} \label {test-fcn5}
\int_{\B^*} | \nabla v|^2 dm & = & n^2 \omega_n^{2/n} \int_0^{|\B^*|}
s^{2- \frac{2}{n}} \left ( \frac{dv^*}{ds} \right )^2 ds \\ \nonumber 
& \leq & \calCp(\Omega) \int_0^{|\B^*|} \left ( - \frac{dv^*}{ds} \right ) 
\int_0^s (v^*(t))^{p-1} dt ds \\ \nonumber 
& = & \calCp(\Omega) \int_0^{|\B^*|}(v^*(t))^{p-1} \int_t^{|\B^*|} 
- \left ( \frac{dv^*}{ds} \right ) ds dt \\ \nonumber 
& = & \calCp(\Omega) \int_0^{|\B^*|} (v^*(t))^p dt. \end {eqnarray}

Combining \eqref{test-fcn4} and \eqref{test-fcn5} we see
\begin {equation} \label {test-fcn6} 
\frac{ \int_{\B^*} |\nabla v|^2 dm}{ \left ( \int_{\B^*} 
|v|^p dm \right )^{2/p}} \leq \calCp(\Omega) \left ( \int_{\B^*}
|v|^p dm \right )^{\frac{p-2}{p}} = \calCp(\B^*) \left ( \int_{\B^*}
|v|^p dm \right )^{\frac{p-2}{p}}.\end {equation}
However, with $p \leq 2$ and $\int_{\B^*} |v|^p dm > 1$ the 
inequality \eqref{test-fcn6} becomes 
$$\frac{\int_{\B^*} |\nabla v|^2 dm} {\left ( \int_{\B^*}
v^p dm \right )^{2/p}} < \calCp(\B^*),$$
which is impossible. 
\end {proof} 

In conclusion, we have shown that there exists $s_1 \in (0,|\B^*|)$ such that 
$\phi^*> u^*$ in $[0,s_1)$ and $u^* \geq \phi^*$ in $(s_1,|\B^*|]$. Extend 
$\phi^*$ to be zero on the interval $(s, |\Omega|]$, so that the inequality 
above continues to hold. That $\phi^*(s) \geq u^*(s)$ for $0 \leq s \leq s_1$ 
immediately implies 
\begin {equation} \label {integral-ineq2} 
\int_0^s (\phi^*)^p d\tau \geq \int_0^s (u^*)^p d\tau \textrm{ for all }
s \in [0, s_1].\end {equation} 
We claim \eqref{integral-ineq2} also holds for $s \in [s_1, |\Omega|]$. To see 
this, let 
$$I(s)  = \int_0^s (\phi^*)^p d\tau - \int_0^s (u^*)^p d\tau,$$
and observe 
$$I(|\Omega|) = I(0) = 0, \qquad I'(s) = (\phi^*(s))^p - (u^*(s))^p.$$
Therefore $I$ is 
nondecreasing in the interval $[0,s_1]$ and nonincreasing in the interval 
$[s_1, |\Omega|]$. In particular, if $I(s_2) < 0$ for some $s_2 \in 
(s_1, |\Omega|)$ then $I(|\Omega|) < 0$ as well, which does not occur. 
Thus \eqref{integral-ineq2} holds for all $s \in [0,|\Omega|]$. 

Combining \eqref{integral-ineq2} with Lemma \ref{HLP-lemma} we 
conclude that for all $q \geq p$ we have 
$$\left ( \int_\Omega u^q dm \right )^{1/q} \leq \left ( 
\int_{\B^*} \phi^q dm\right )^{1/q} \leq K \left ( \int_{\B^*} \phi^p dm 
\right )^{1/p} = K \left ( \int_\Omega u^p dm \right )^{1/p},$$
where 
$$K = \frac{ \left ( \int_{\B^*} \phi^p dm \right )^{1/p}} 
{\left ( \int_{\B^*} \phi^q dm \right )^{1/q}}$$
and we have used the normalization 
$$\int_{\B^*} \phi^p dm = \int_\Omega u^p dm . $$
\hfill $\square$


\begin {thebibliography}{999}


\bibitem{vdB} M. van den Berg, \textsl{Estimates for
the torsion function and Sobolev constants.\/} Potential Analysis, {\bf 36} (2012) 607--616. 

\bibitem {vdBC} M. van den Berg and T. Carroll. \textsl{Hardy inequality
and $L^p$ estimates for the torsion function.\/} Bull. London Math. Soc. 
(2009), 980--986. 



\bibitem{CRpPR1}T.\ Carroll and J.\ Ratzkin,
\textsl{Two isoperimetric inequalities for the Sobolev constant.\/}
Z.\ Angew.\ Math.\ Phys.\ {\bf 63} (2012), 855--863. 

\bibitem{CRpPR2} T.\ Carroll and J.\ Ratzkin, 
\textsl{An isoperimetric inequality for extremal Sobolev functions.\/}
RIMS Kokyuroku Bessatsu {\bf B43} (2013), 1--16. 


\bibitem {Chiti1}G.\ Chiti, \textsl{An isomerimetric inequality for the 
eigenfunctions of linear second order elliptic equations.\/} 
Boll.\ Un.\ Mat.\ Ital.\ A {\bf 1} (1982), 145--151. 

\bibitem{Chiti2}G.\ Chiti,
\textsl{A reverse H\"older inequality for the eigenfunctions of 
linear second order elliptic operators.\/}
Z.\ Angew.\ Math.\ Phys.\ {\bf 33} (1982), 143--148.

\bibitem {GNN} B.\ Gidas, W.\-N.\ Ni, and L.\ Nireberg, \textsl{Symmetry 
and related properties via the maximum principle}. Comm. Math. Phys. 
{\bf 68} (1979), 209--243. 

\bibitem {HLP} G.\ H.\ Hardy, J.\ E.\ Littlewood, and G.\ P\' olya. \textsl{Some 
simple inequalities satisfied by convex functions.} Messenger Math. 
{\bf 58} (1929), 145--152. 

\bibitem{KJ}M.-T. Kohler-Jobin,
\textsl{Sur la premi\`ere fonction propre d'une membrane: une extension 
\`a $N$ dimensions de l'in\'egalit\'e iso\-p\'eri\-m\'e\-trique de Payne-Rayner.\/}
Z.\ Angew.\ Math.\ Phys.\ {\bf 28} (1977), 1137--1140.

\bibitem{KJ3}M.-T. Kohler-Jobin,
\textsl{Isoperimetric monotonicity and isoperimetric inequalities of 
Payne-Rayner type for the first eigenfunction of the Helmholtz problem.\/}
Z.\ Angew.\ Math.\ Phys.\ {\bf 32} (1981), 625--646. 

\bibitem{PR} L.\ Payne and M.\ Rayner, 
\textsl{An isoperimetric inequality for the first 
eigenfunction in the fixed membrane problem.} 
Z.\ Angew.\ Math.\ Phys.\ {\bf 23} (1972), 13--15.

\bibitem{PR2}L.\ Payne and M.\ Rayner, 
\textsl{Some isoperimetric norm bounds for solutions of the Helmholtz 
equation.\/} Z.\ Angew.\ Math.\ Phys.\ {\bf 24} (1973), 105--110.

\bibitem{PS} G. P\' olya and G. Szeg\H o,
{\em Isoperimetric Inequalities in Mathematical Physics}. 
Princeton University Press (1951).

\bibitem {Tal} G. Talenti, \textsl{Elliptic equations and rearrangements.\/}
Ann.\ Scuola\ Norm.\ Sup.\ Pisa\ Cl.\ Sci.\ {\bf 3} (1976), 697--718. 

\end {thebibliography}
\end{document}